# Models Parametric Analysis
# via Adaptive Kernel Learning

Vladimir Norkin [0000-0003-3255-0405] and Alois Pichler [0000-0001-8876-2429]

V.M. Glushkov Institute of Cybernetics
of the National Academy of Sciences of Ukraine, Kyiv, 03178 Ukraine &
National Technical University of Ukraine "Igor Sikorsky Kyiv Polytechnic Institute"

vladimir.norkin@gmail.com

Chemnitz University of Technology, Faculty of Mathematics, Chemnitz, Germany

alois.pichler@math.tu-chemnitz.de

**Abstract.** Any applied mathematical model contains parameters. The paper proposes to use kernel learning for the parametric analysis of the model. The approach consists in setting a distribution on the parameter space, obtaining a finite training sample from this distribution, solving the problem for each parameter value from this sample, and constructing a kernel approximation of the parametric dependence on the entire set of parameter values. The kernel approximation is obtained by minimizing the approximation error on the training sample and adjusting kernel parameters (width) on the same or another independent sample of parameters. This approach to learning complex dependencies is called kernel learning (or kernel SVM). Traditionally, the kernel learning is considered in the so-called Reproducing Kernel Hilbert Space (RKHS) with a fixed kernel. The novelty of our approach is that we consider the kernel learning in a broad subspace of square integrable functions with a corresponding the $L_2$-norm regularization. This subspace contains linear combinations of kernel functions with kernels of different shape at different data points. The approach essentially uses a derived analytical representation of the $L_2$-norm for kernel functions. Thus the approach substantially extends the flexibility of the traditional kernel SVM for account of the number of adjusted parameters and minimization of the training error not only over weights of kernels but also over their shapes. The important issue of selection of the optimal regularization parameter is resolved by minimizing the test error over this parameter. Numerical illustrations are provided.

# 1    Introduction

Machine learning methods are gradually penetrating into all areas of applied mathematics. They allow one to approximate complex non-linear mathematical relationships by learning from examples. Any applied mathematical problem contains parameters. In practice, it is desirable to solve this problem not only for a specific value of the parameter vector, but also for the entire spectrum of possible parameter values. This is necessary both to study the dependence of the solution of the problem on the parameters, and to understand the stability of the solution in the case of variation (perturbation) of the parameters and to understand the risks due to possible inaccurate knowledge (uncertainty) in the values of the parameters.

Traditionally, parametric analysis of a mathematical model consists in studying the stability, sensitivity and robustness of the solution with respect to changing the parameters of the problem. For this purpose, countless special methods and approaches have been developed to study the solution dependence on certain parameters. However, the problem is that there can be a lot of parameters (e.g., a matrix), and they coherently affect the solution nonlinearly.

An important application of kernel method is visualization of spatial data [1, 2], where the selection of kernel width is essential for smooth representation of discrete spatial data.

Another classical approach to the approximation of an unknown dependence is being developed in multivariate statistics in the form of parametric and nonparametric regression analysis, where the main emphasis is on the objective statistical nature of the initial data [3-6].

This paper proposes to use machine learning methods to approximate the multidimensional nonlinear dependence of the solution of the original problem on the parameters of the problem. The idea of the approach is to set a (uniform or some subjective) distribution on the problem parameter space, obtain a finite training sample from this distribution, solve the original problem for each parameter value from this sample, and construct a kernel or neural network approximation of the parametric dependence on the entire set of parameter values. In this case, the approximation coefficients are chosen by minimizing the approximation error on the sample and checking the accuracy of the obtained approximation on another independent test sample of parameters. This approach can be implemented both within the framework of the kernel Support Vector Machine (kernel SVM) [7-9] and with the help of deep neural networks (DNN) [10].

To implement this approach, it is necessary to adapt machine learning methods to the problem of learning parametric dependencies. They must be flexible enough to be able to approximate complex multivariate non-smooth and possibly discontinuous dependencies.

Adaptation of kernels to data has long history and is an active field of research in approximation theory, nonparametric statistics and machine learning [11-22]. A review of kernel based methods in different fields can be found in [20]. In particular, paper [17]  introduces Local-Adaptive-Bandwidth RBF kernels (LAB RBF), which are similar to the anisotropic kernels, proposed in the present paper but from a differ-ent function approximation (smoothing) approach. Paper [19] explores application of global optimization methods for finding optimal univariate kernel shape parameter within the leave-one out cross validation framework [22]. Papers [21, 23, 24] extend theory of Support Vector Machines in Reproducing Kernel Hilbert Spaces (RKHS) to SVM in Reproducing Kernel Banach Spaces (RKBS).

The present paper, first extends the kernel support vector machine by using adaptive kernels similar to [17], i.e. kernels with spatially variable widths or non-uniformly anisotropic kernels. Moreover, the parameters of the kernels used in the approximation can be different and depend on the sample. The number of adjustable parameters becomes much larger than the sample size, which makes the kernel support vector machine approaching in flexibility to multi-layer neural networks. However, a new learning problem arises for such an adaptive support vector machine.

Second, the paper proposes to consider a regression problem not in RKHS [25, 26] but in a kernel subspace of $L_2(\mathbb{R}^n)$ with a standard inner product in this space and with a standard $L_2$-norm in the regularization term. Unlike [17], where the kernel matrix is obtained through the so-called kernel trick, we derive kernel matrices by analytical calculation of $L_2$-norm of kernel functions.

Third, selection of the regularization parameter is done through minimizing absolute error on training and test data over this parameter.

Thus the contribution of the paper includes the following.

The paper proposes and explores a wide area of application of machine learning methods in applied mathematics, namely, parametric analysis of mathematical models by machine learning methods, in particular, by adaptive kernel support vector machine. Thus the paper makes a bridge between classical mathematical modeling and new machine learning techniques.

The kernel support vector machine is extended by the ability to use a continuum (multivariate parametric) family of kernels to approximate the dependence under study. For example, it is proposed to use variable width kernels and variable anisotropic kernels to reduce the approximation error. Thus the standard for kernel leaning Reproducing Kernel Hilbert Space is replaced by much broader kernel subspace of square integrable functions. While kernels in RKHS act on functions (via the inner product) as Dirac's delta-function, kernels in $L_2$ act on functions as smoothing/mollifying operators [27].

An essential problem in machine learning is selection of regularization, in particular the scale of the regularization parameter. In the paper we select the regularization parameter by minimization of the approximation error on training and/or test data.

To train the adaptive support vector machine, it is proposed to use a combination of traditional kernel weights optimization methods (quadratic as in SVM) and global optimization methods to optimize kernel parameters.

## 2 Nonparametric analysis of parametric models

Any applied mathematical model contains parameters. The paper proposes to use machine learning methods based on kernel functions for parametric model analysis. The approach consists in specifying some distribution on the parameter space, obtaining a finite training parameter sample from this distribution, solving the model for each selected value of the parameters, and constructing a kernel approximation of the parametric dependence for the entire set of parameters. The kernel approximation is obtained by minimizing the approximation error on the training sample and testing the approximation on another independent test sample of the parameters. The novelty of the approach lies in the adaptation of the shape of the kernels to the observations. This makes the flexibility of the kernel method comparable to the flexibility of neural networks. In addition, the training set of parameters (and corresponding solutions) can be successively expanded, and the corresponding approximation of the parametric dependence can be refined online using adaptive stochastic quasi-gradient methods.

EXAMPLE (a toy parametric problem). Consider an input-output type model in the form of a quadratic equation $y = x^2 + ax + b$, where $x$ is the input of the model, and $y$ is the output of the model. Suppose we are looking for an input $x$, at which the output is zero, i.e. $y = 0$. For this, we need to solve the equation $x^2 + ax + b = 0$. Its solution $x^* = x^*(a,b) = f(a,b)$ depends on the model parameters $a \in [-2, 2]$, $b \in [-2, 2]$. The solution can be found analytically,

$$f^{\pm}(a,b) = x^{\pm}(a,b) = \frac{-a \pm \sqrt{a^2 - 4b}}{2}, \tag{1}$$

and it can also be approximated nonparametrically, e.g., statistically. Namely, for a random sample $(a_1, b_1), ..., (a_m, b_m)$, we find the appropriate solutions $x^{\pm}(a_1, b_1), ..., x^{\pm}(a_m, b_m)$, and for arbitrary parameters $(a, b)$ we build a regression dependence $x^{\pm}(a, b)$ based on the data $x_1^{\pm} = x^{\pm}(a_1, b_1), ..., x_N^{\pm} = x^{\pm}(a_m, b_m)$ using the SVM method. Then we compare the obtained approximation with the exact solution.

### 2.1 Nadaraya-Watson kernel regression

One option to reconstruct general parametric dependences is provided by the Nadaraya-Watson method [4-6]. Suppose there is a set of observations $\{\overline{y}_i\}_{i=1}^m$ of the unknown functional dependence $y = F(x) : \mathbb{R}^n \to \mathbb{R}$ at points $\{\overline{x}_i\}_{i=1}^m$. The problem is to restore function $y = F(x)$ on the basis of observations $\{\overline{x}_i, \overline{y}_i\}_{i=1}^m$. Nadaraya-Watson kernel estimate (NW-regression) of function $F(\cdot)$ is found as a solution of the problem [4-6]:

$$\Phi_{\theta}^m(x,z) = \frac{1}{2N} \sum_{i=1}^m \frac{1}{\sigma^r} K\left(\frac{\overline{x}_i - x}{\sigma}\right) \|z - \overline{y}_i\|^2 \xrightarrow[z \in \mathbb{R}^s]{} \min,$$

where $k_\sigma(\bar{x}_i, x) = \dfrac{1}{\sigma^r} K\left(\dfrac{\bar{x}_i - x}{\sigma}\right)$ is some density function [21], $\sigma > 0$ is the kernel (width) parameter. The solution of this problem is

$$z_\sigma(x) = \sum_{i=1}^m \bar{y}_i K\left(\frac{\bar{x}_i - x}{\sigma}\right) \bigg/ \sum_{i=1}^m K\left(\frac{\bar{x}_i - x}{\sigma}\right). \qquad (2)$$

The number of essential summands in (2) for a given $x$ depends on the kernel width $\sigma$ and on the shape of the kernel. One can explicitly regulate the number of summands in the $k$-nearest neighbor Nadaraya-Watson ($k$NNW) regression by finding the set $I_k(x)$ of $k$ nearest to $x$ points of the set $\{\bar{x}_i\}_{i=1}^N$ and calculating

$$z_{\sigma,k}(x) = \sum_{i \in I_k(x)} \bar{y}_i K\left(\frac{\bar{x}_i - x}{\sigma}\right) \bigg/ \sum_{i \in I_k(x)} K\left(\frac{\bar{x}_i - x}{\sigma}\right). \qquad (3)$$

If $k = m$, then (3) coincides with (2). However, in this case function (3) can be discontinuous even for a continuous kernel $K(\cdot)$.

Nadaraya-Watson regression is rather robust to the choice of kernel width $\sigma$, since $z_\sigma(x)$ is just some average value of quantities $\{\bar{y}_i\}_{i=1}^m$. There is a room for further modification of the NW-regression by adjusting parameter $\sigma$ to dada. For example, if $K(\cdot)$ is an one-dimensional function, $K\left(\dfrac{\bar{x}_i - x}{\sigma}\right) = K\left(\dfrac{\|\bar{x}_i - x\|}{\sigma}\right)$ and $\sigma = \sigma(\bar{x}_i, x) = \|\bar{x}_i - x\|$, then (3) becomes just an arithmetic average of $\bar{y}_i$, $i \in I_k(x)$,

$$z_k(x) = \frac{1}{k} \sum_{i \in I_k(x)} \bar{y}_i. \qquad (4)$$

Remark also that Nadaraya-Watson method is applicable to vector regression $y = F(x): \mathbb{R}^n \to \mathbb{R}^l$, when $\bar{y}_i \in \mathbb{R}^l$, $l \ge 1$, and to the $l$-class classification problem, when $\bar{y}_i \in \{0,1\}^l$ are Boolean vectors or $\bar{y}_i \in [0,1]^l$ are probability distributions with the sum of components equal to one.

EXAMPLE (continued). The kernel, for example, may have the form

$$K\left(\frac{a - \bar{a}}{\sigma_a}, \frac{b - \bar{b}}{\sigma_b}\right) = \exp\left(-\frac{|a - \bar{a}|^2}{\sigma_a^2} - \frac{|b - \bar{b}|^2}{\sigma_b^2}\right) = \exp\left(-\frac{|a - \bar{a}|^2}{\sigma_a^2}\right) \times \exp\left(-\frac{|b - \bar{b}|^2}{\sigma_b^2}\right),$$

then Nadaraya-Watson kernel regression is

$$z_m^\pm\left(a, b; \sigma_a, \sigma_b; \{(a_i, b_i)\}_{i=1}^m\right) = \sum_{i=1}^m x_i^\pm \frac{K\left(\dfrac{a - a_i}{\sigma_a}, \dfrac{b - b_i}{\sigma_b}\right)}{\sum_{i'=1}^m K\left(\dfrac{a - a_{i'}}{\sigma_a}, \dfrac{b - b_{i'}}{\sigma_b}\right)}.$$

Now it is necessary to compare functions $f^{\pm}(a,b)$ and $z_m^{\pm}\left(a,b;\sigma_a,\sigma_b;\{(a_i,b_i)\}_{i=1}^m\right)$, i.e. to calculate the approximation error

$$\Delta_m^{\pm}\left(\sigma_a,\sigma_b;\{(a_i,b_i)\}_{i=1}^m\right)=\sup_{(a,b)\in[-2,2]^2}\left|f^{\pm}(a,b)-z_n^{\pm}\left(a,b;\sigma_a,\sigma_b;\{(a_i,b_i)\}_{i=1}^m\right)\right|.$$

Here, the kernel width parameters $\sigma=(\sigma_a,\sigma_b)$ are assumed to be fixed and preselected. Is it possible to improve (reduce) the error at the expense of variables $(\sigma_a,\sigma_b)$, i.e. at the expense of adapting the kernels (their width) to the data? For example, let $(\sigma_a,\sigma_b)$ depend on $i$, that is, consider $\left\{(\sigma_a^i,\sigma_b^i)\right\}_{i=1}^m$. But how to choose these values?

The generalized Nadaraya-Watson regression function now looks like:

$$z_m^{\pm}\left(a,b;\{\sigma_a^i,\sigma_b^i\}_{i=1}^m;\{(a_i,b_i)\}_{i=1}^m\right)=\sum_{i=1}^m x_i^{\pm}\frac{K\left(\dfrac{a-a_i}{\sigma_a^i},\dfrac{b-b_i}{\sigma_b^i}\right)}{\sum_{i'=1}^m K\left(\dfrac{a-a_{i'}}{\sigma_a^{i'}},\dfrac{b-b_{i'}}{\sigma_b^{i'}}\right)}.$$

Suppose we have an additional data set $\left\{x_j^{\pm},(a_j,b_j)\right\}_{j=m+1}^{m'}$. Now we will find the better kernel width parameters $\left\{(\sigma_a^i,\sigma_b^i)\right\}_{i=1}^m$ by minimizing the additional error function

$$\Delta_m^{m'}\left(\{(\sigma_a^i,\sigma_b^i)\}_{i=1}^m\right)=\sum_{j=m+1}^{m'}\left(x_j^{\pm}-z_m^{\pm}\left(a_j,b_j;\{(\sigma_a^i,\sigma_b^i)\}_{i=1}^m;\{(a_i,b_i)\}_{i=1}^m\right)\right)^2\to\min_{\{(\sigma_a^i,\sigma_b^i)\}_{i=1}^m}.$$

This is a high-dimensional optimization problem with respect to variables $\left\{(\sigma_a^i,\sigma_b^i)\right\}_{i=1}^m$. It is non-convex and (possibly) has many extremes. Various optimization methods can be used to solve it, including gradient methods and stochastic gradient methods.

Thus, the minimization of the error $\Delta_m^{m'}\left(\{(\sigma_a^i,\sigma_b^i)\}_{i=1}^m\right)$ by variables $\left\{(\sigma_a^i,\sigma_b^i)\right\}_{i=1}^m$ should be done using the same methods as for the training neural networks with the use of automatic differentiation programs and means of parallelizing calculations.

By using kernels with variable width, we expect to obtain better approximations of the parametric solutions of our model or to reduce the training sample.

## 2.2 The kernel SVM for a regression problem

Standard regression problem in a Reproducing Kernel Hilbert Space (RKHS) with kernel $k(\cdot,\cdot)$ has the following form [8, 28]. Let $\left\{(\bar{x}_i,\bar{y}_i),\,i=1,...,m\right\}$ be a sample from some (unknown) distribution $\rho(x,y)$. To find the conditional mean function $F(x) = \int z\rho(x,y)P(dy) \big/ \int \rho(x,y)P(dy)$ one considers the so-called empirical risk minimization problem (kernel SVM)

$$R^m(f) = \sum_{i=1}^{m}\left(\overline{y}_i - f(\overline{x}_i)\right)^2 + \lambda\left\|f(\cdot)\right\|_{H_k}^2 \to \min_{f\in H_k}, \qquad (5)$$

where $f(\cdot)$ is an approximate searched function, $H_k$ is a corresponding Reproducing Kernel Hilbert Space (RKHS) [8, 25] with reproducing kernel $k(x,\overline{x})$, $\lambda$ is a regularization parameter, $\left\|f\right\|_{H_k}$ is the norm of function $f \in H_k$. It is known that a solution for this optimization problem has the form [8] $f(\cdot) = \sum_{i=1}^{m}\alpha_i k(\overline{x}_i, \cdot)$ and its norm is

$$\left\|f\right\|_{H_k} = \left(\sum_{i,j=1}^{m}\alpha_i\alpha_j k(\overline{x}_i,\overline{x}_j)\right)^{1/2} = \left(\sum_{i,j=1}^{m}\alpha_i\alpha_j k_{ij}\right)^{1/2}$$

with some real coefficients $\alpha = \left(\alpha_1,...,\alpha_N\right)^T$, $k_{ij} = k(\overline{x}_i,\overline{x}_j)$. Thus the regression problem (5) is reduced to the finite dimensional optimization problem

$$R^m(\alpha) = \sum_{s=1}^{m}\left(y_s - \sum_{t=1}^{m}\alpha_t k_{st}\right)^2 + \lambda\sum_{s,t=1}^{m}\alpha_s\alpha_t k_{st} \to \min_{\alpha}. \qquad (6)$$

This is a quadratic unconstrained optimization problem. Its solution can be found in a closed form [17, 29] $\alpha = \left(K + \lambda I\right)^{-1}\overline{y}$, where matrix $K = \left\{k(\overline{x}_i,\overline{x}_j)\right\}_{i,j=1}^{m}$ is assumed symmetric and positive definite, $I$ is the identity matrix, $\left(K + \lambda I\right)^{-1}$ is the inverse of matrix $K + \lambda I$, $\overline{y} = \left(\overline{y}_1,...,\overline{y}_m\right)^T$. The optimal value for problem (6) is $R_{\min}^m = \lambda\overline{y}^T\left(K + \lambda I\right)^{-1}\overline{y}$ [17, 29].

**Remark 1.** In [29] the kernel estimate of the regression function had the form $f(\cdot) = m^{-1}\sum_{i=1}^{m}\alpha_i k(\overline{x}_i, \cdot)$ and the error minimization problem had the form

$$\frac{1}{m}\sum_{i=1}^{m}\left(\overline{y}_i - f(\overline{x}_i)\right)^2 + \lambda\left\|f(\cdot)\right\|_{H_k}^2 \to \min_{f\in H_k}$$

so the optimal weights and the optimal value for (6) were $\alpha = \left(m^{-1}K + \lambda I\right)^{-1}\overline{y}$ and $R_{\min}^m = \lambda m^{-1}\overline{y}^T\left(m^{-1}K + \lambda I\right)^{-1}\overline{y}$.

The standard kernel learning method assumes using kernels with fixed (width) parameters although dependent on the size $m$ of the sample $\left\{(\overline{x}_i,\overline{y}_i)\right\}_{i=1}^{m}$. In recursive kernel learning, see, e.g. [30] (and references there in), the kernel width is changed with the course of iterations and thus is different for kernels sitting at different sample points. In the present paper we further adapt kernels to data by admitting different width of kernels at different points.

Remark that nonquatratic loss functions can be used in (7). Then SVM-estimate is solution of the following infinite dimensional optimization problem

$$\sum_{i=1}^{m} c\left(\overline{y}_i, f(\overline{x}_i)\right) + \lambda \left\| f(\cdot) \right\|_{H_k}^2 =$$

$$= \sum_{i=1}^{m} c\left(\overline{y}_i, \sum_{j=1}^{m} \alpha_j k(\overline{x}_j, \overline{x}_i)\right) + \lambda \sum_{i,j=1}^{m} \alpha_i \alpha_j k(\overline{x}_j, \overline{x}_i) \to \min_{\alpha},$$

where $c(\cdot, \cdot)$ is some loss function, e.g. $c(y, y') = |y - y'|^2$ or $c(y, y') = |y - y'|$ [8, 28], $\lambda > 0$ is a regularization parameter, $\left\| f(\cdot) \right\|_{H_k}^2 = \sum_{i,j=1}^{m} \alpha_i \alpha_j k(\overline{x}_i, \overline{x}_j)$.

# 3 Exposition of the novel solution approach

In this section we extend the standard kernel SVM in three aspects:

First, we consider the regression problem not in a RKHS but in $L_2(\mathbb{R}^n)$;

Second, we develop an adaptive kernel SVM, where kernels are adapted to data;

Third, we suggest choosing the regularization parameter through cross validation technique.

## 3.1 Multidimentional kernel regression in $L_2(\mathbb{R}^n)$

Let $L_2(\mathbb{R}^n)$ be a space of square integrable functions on $\mathbb{R}^n$. For two functions in this space, $f(\cdot)$ and $g(\cdot)$, the $L_2(\mathbb{R}^n)$-inner product is defined as

$$\langle f, g \rangle = \int_{\mathbb{R}^n} f(x) g(x) dx.$$

Let $\rho(x) \geq 0, x \in \mathbb{R}^n$, be some even probabilistic density function, $\int_{\mathbb{R}^n} \rho(x) dx = 1$, $\rho(x) = \rho(-x)$. Consider a family of kernels $k_\sigma(x, \overline{x}) = \sigma^{-n} \rho\left(\sigma^{-1}(x - \overline{x})\right)$, $\sigma > 0$, and the inner products $\left\langle k_{\sigma_1}(x, \overline{x}_1), k_{\sigma_2}(x, \overline{x}_2) \right\rangle = \int_{\mathbb{R}} k_{\sigma_1}(x, \overline{x}_1) k_{\sigma_2}(x, \overline{x}_2) dx$.

Let function $f(\cdot) \in L_2(\mathbb{R}^n)$ and $f(x) = 0$ for $x \notin X$, $X$ is a bounded Borell measurable subset in $\mathbb{R}^n$. Define the averaged function

$$f_\sigma(x) = \int_{\mathbb{R}^n} f(y) k_\sigma(y - x) dy = \int_{\mathbb{R}^n} f(x + \sigma z) \rho(z) dz.$$

If $f(\cdot)$ is continuous in a vicinity of $x$, then $\lim_{\sigma \to +0} f_\sigma(x) = f(x)$, i.e. functions $f_\sigma(\cdot)$ approximate $f(\cdot)$. Moreover, if $f(x)$ is Lipschitz continuous in $x \in \mathbb{R}^n$ with constant $L$ and $\int_{\mathbb{R}^n} \|z\| \rho(z) dz = C < +\infty$, then

$$\left| f(x) - f_\sigma(x) \right| \leq \int_{\mathbb{R}^n} \left| f(x) - f(x + \sigma z) \right| \rho(z) dz \leq \sigma L C .$$

Let $\left\{ \overline{x}_i \right\}_{i=1}^m$ be a sample of points uniformly distributed in $X$, denote $\overline{y}_i = f(\overline{x}_i)$. Then the Monte Carlo estimate $\overline{f}_\sigma^m(x)$ of $f_\sigma(x)$ has the form

$$\overline{f}_\sigma^m(x) = m^{-1} \sum_{i=1}^m \overline{y}_i \, k_\sigma(\overline{x}_i - x) .$$

Thus for sufficiently large $m$ and sufficiently small $\sigma = \sigma(m)$ function $f(x)$ can be approximated by the kernel estimates

$$\overline{f}_\sigma^m(x) = \sum_{i=1}^m \alpha_i k_\sigma(\overline{x}_i - x), \qquad \alpha_i = \overline{y}_i / m . \tag{8}$$

The convergence of $\overline{f}_{\sigma_m}^m(x)$ will be the subject of separate investigation but in the present paper we consider the following question:

Can the estimate $\overline{f}_\sigma^m(x)$ be improved for account of choice of coefficient $\alpha_i$ other than $\alpha_i = \overline{y}_i / m$ and for account of choice of $\sigma$ dependent on the sample $\left\{ (\overline{x}_i, \overline{y}_i) \right\}_{i=1}^m$?

We partially answer this question by conducting numerical experiments with problem considered in the Example of Section 2. First we give a method how to select coefficients $\alpha_i$ in (8). We develop this method by analogy with kernel SVM regression in RKHS.

Consider the subspace $KL_2$ of kernel functions of the form

$$f_{\alpha, \overline{x}, \sigma}(x) = \sum_{i=1}^m \alpha_i k_{\sigma_i}(\overline{x}_i, x) , \quad \overline{x} = \left\{ \overline{x}_i \in \mathbb{R}^n, \, i = 1, ..., m \right\},$$

$$\alpha = \left( \alpha_1, ..., \alpha_m \right)^{\mathrm{T}} \in \mathbb{R}^m, \quad \sigma = \left( \sigma_1, ..., \sigma_m \right) \in S^m, \quad m \geq 1 .$$

If functions $f(\cdot)$ and $g(\cdot)$ have the forms

$$f(\cdot) = \sum_{i=1}^m \alpha_i k_{\sigma_i}(\overline{x}_i, \cdot) , \qquad g(\cdot) = \sum_{j=1}^l \beta_j k_{\mu_j}(\hat{x}_j, \cdot) ,$$

then the corresponding inner product takes on the form

$$\begin{aligned}
\langle f(\cdot), g(\cdot) \rangle &= \left\langle \sum_{i=1}^m \alpha_i k_{\sigma_i}(\overline{x}_i, \cdot), \sum_{j=1}^l \beta_j k_{\mu_j}(\hat{x}_j, \cdot) \right\rangle \\
&= \sum_{i=1}^m \sum_{j=1}^l \alpha_i \beta_j \left\langle k_{\sigma_i}(\overline{x}_i, \cdot), k_{\mu_j}(\hat{x}_j, \cdot) \right\rangle \\
&= \sum_{i=1}^m \sum_{j=1}^l \alpha_i \beta_j \int_{\mathbb{R}^n} k_{\sigma_i}(\overline{x}_i, x) k_{\mu_j}(\hat{x}_j, x) dx \\
&= \sum_{i,j=1}^{m,l} \alpha_i \beta_j k_{ij} = \alpha^{\mathrm{T}} K \beta ,
\end{aligned}$$

where $K = \left\{ k_{ij} \right\}_{i,j=1}^{m,l}$ is a matrix (to be more specified later) with entries

$$k_{ij} = \left\langle k_{\sigma_i}(\overline{x}_i,\cdot), k_{\mu_j}(\hat{x}_j,\cdot) \right\rangle = \int_{\mathbb{R}^n} k_{\sigma_i}(\overline{x}_i,x) k_{\mu_j}(\hat{x}_j,x) dx .$$

The norm of $f_{\alpha,\overline{x},\sigma}(x)$ in this subspace is:

$$\left\| f_{\alpha,\overline{x},\sigma}(\cdot) \right\|_{L_2} = \left( \sum_{i,j=1}^{m} \alpha_i \alpha_j \left\langle k_{\sigma_i}(\overline{x}_i,\cdot), k_{\sigma_j}(\overline{x}_j,\cdot) \right\rangle \right)^{1/2} = \left( \alpha^\mathsf{T} K \alpha \right)^{1/2} ,$$

where vector $\alpha = (\alpha_1, \alpha_2, ...)^\mathsf{T} \in \mathbb{R}^m$, and the matrix

$$K = \begin{pmatrix} \left\langle k_{\sigma_1}(\overline{x}_1,\cdot), k_{\sigma_1}(\overline{x}_1,\cdot) \right\rangle & ... & \left\langle k_{\sigma_1}(\overline{x}_1,\cdot), k_{\sigma_m}(\overline{x}_m,\cdot) \right\rangle \\ ... & \left\langle k_{\sigma_j}(\overline{x}_j,\cdot), k_{\sigma_j}(\overline{x}_j,\cdot) \right\rangle & ... \\ \left\langle k_{\sigma_m}(\overline{x}_m,\cdot), k_{\sigma_1}(\overline{x}_1,\cdot) \right\rangle & ... & \left\langle k_{\sigma_m}(\overline{x}_m,\cdot), k_{\sigma_m}(\overline{x}_m,\cdot) \right\rangle \end{pmatrix}. \quad (9)$$

Since matrix $K$ is symmetric and the coefficients on the main diagonal are positive,

$$\left\langle k_{\sigma_j}(\overline{x}_j,\cdot), k_{\sigma_j}(\overline{x}_j,\cdot) \right\rangle = \int_{\mathbb{R}} k_{\sigma_j}(\overline{x}_j,\cdot) k_{\sigma_j}(\overline{x}_j,\cdot) dx = \frac{1}{\sigma_j^{2n}} \int_{\mathbb{R}^n} \rho^2 \left( \frac{\overline{x}_j - x}{\sigma_j} \right) dx > 0 ,$$

then $K$ is positive definite but may be ill-posed.

The idea (and the novelty) of the our approach consists in the extension of the standard functional space RKHS to a kernel subspace of $L_2(\mathbb{R}^n)$, containing functions

$$f_{\alpha,\overline{x},\sigma}(x) = \sum_{i=1}^{m} \alpha_i k_{\sigma_i}(\overline{x}_i,x) , \quad \overline{x} = (\overline{x}_1,...,\overline{x}_m)^\mathsf{T} \in \mathbb{R}^m ,$$

$$\alpha = (\alpha_1,...,\alpha_m)^\mathsf{T} \in \mathbb{R}^m , \quad \sigma = \{\sigma_1,...,\sigma_m\} , \quad m \geq 1 ,$$

with kernel functions $k_{\sigma_i}(\overline{x}_i,x)$ specific for each sample point $\overline{x}_i$. This space is much larger than RKHS with a fixed kernel.

The kernel function $k_{\sigma}(x,\overline{x})$, $x,\overline{x} \in \mathbb{R}^n$, can be, e.g., the Gaussian kernel

$$k_{\sigma}(x,\overline{x}) = \frac{1}{(2\pi)^{n/2} |\Sigma|^{1/2}} e^{-\frac{1}{2}(x-\overline{x})^\mathsf{T} \Sigma^{-1}(x-\overline{x})} ,$$

where $\Sigma$ is a positive definite matrix, $|\Sigma|$ is the determinant of $\Sigma$. As particular cases one can take $\Sigma = I\sigma^2$, $\sigma > 0$, or $\Sigma = I \times (\sigma_1^2,...,\sigma_n^2)^\mathsf{T}$, where $I$ is the identity diagonal matrix. The corresponding kernel regression functions have the form:

$$f_{\alpha,\overline{x},\sigma}(x) = \sum_{i=1}^{m} \frac{\alpha_i}{(2\pi)^{n/2} \sigma_i^n} e^{-\|x-\overline{x}_i\|^2/(2\sigma_i^2)} \qquad \text{(with variable kernel widths}$$

$\sigma = (\sigma_1,...,\sigma_m)$ ),

$$f_{\alpha,\overline{x},\sigma}(x) = \sum_{i=1}^{m} \frac{\alpha_i}{(2\pi)^{n/2}\sigma_{i1}\times...\times\sigma_{in}} e^{-\sum_{l=1}^{n}((x)_l-(\overline{x}_i)_l)^2/(2\sigma_{il}^2)}$$ (with anisotropic variable kernel widths $\sigma = \left(\left\{\sigma_{1j}\right\}_{j=1}^{n},...,\left\{\sigma_{mj}\right\}_{j=1}^{n}\right)$.

The corresponding kernel regression problem for a sample $\left\{(\overline{x}_i, \overline{y}_i), i=1,...,m\right\}$ takes on the form

$$R^m(\alpha,\sigma) = \sum_{i=1}^{m} c\left(\overline{y}_i, f_{\alpha,\overline{x},\sigma}(x_i)\right) + \lambda\Psi\left(f_{\alpha,\overline{x},\sigma}(\cdot)\right) =$$
$$= \sum_{i=1}^{m} c\left(\overline{y}_i, \sum_{j=1}^{m} \alpha_j k_{\sigma_j}(\overline{x}_i, \overline{x}_j)\right) + \lambda\Psi\left(f_{\alpha,\overline{x},\sigma}(\cdot)\right) \to \min_{\alpha,\sigma}, \quad (10)$$

where $\Psi\left(f_{\alpha,\overline{x},\sigma}(\cdot)\right)$ is some regularization term (to be defined, e.g. $\Psi\left(f_{\alpha,\overline{x},\sigma}(\cdot)\right) = \left\|f_{\alpha,\overline{x},\sigma}(\cdot)\right\|_{L_2}^2$). This optimization problem has much more optimization variables than $m$, e.g. in the isotropic case it has $2m$ variables, $\alpha = (\alpha_1,...,\alpha_m)^{\mathrm{T}}$ and $\sigma = (\sigma_1,...,\sigma_m)$. This problem is usually convex in variables $\alpha$ and nonlinear in variables $\sigma$. The optimization procedure for (10) may be consisting of the repeated solutions of the problems

$$R^m(\alpha^*,\sigma) = \min_{\alpha}\left[R^m(\alpha,\sigma) = \left(\sum_{i=1}^{m} c\left(\overline{y}_i, f_{\alpha,\overline{x},\sigma}(\overline{x}_i)\right) + \lambda\Psi\left(f_{\alpha,\overline{x},\sigma}(\cdot)\right)\right)\right], \quad (11)$$

$$R^m(\alpha^*,\sigma^*) = \min_{\sigma}\left[R^m(\alpha^*,\sigma) = \left(\sum_{i=1}^{m'} c\left(\overline{y}_i', f_{\alpha^*,\overline{x},\sigma}(\overline{x}_i')\right) + \lambda\Psi\left(f_{\alpha^*,\overline{x},\sigma}(\cdot)\right)\right)\right], \quad (12)$$

where $\left\{(\overline{x}_i', \overline{y}_i'), i=1,...,m'\right\}$ is some additional sample of observations of function $F$. Here optimization over variables $\alpha \in \mathbb{R}^m$ can be made by some standard software and the (global) optimization over $\sigma$ could be done by deterministic or stochastic adaptive gradient methods [31, 32] and by their nonsmooth versions like this is done in deep learning [33]. A similar procedure for solving problem (10) is applied in [17].

### 3.2 Analytical calculation of the norm

Let us consider a particular case of the kernel regression problem (10) with a quadratic regularization term:

$$R^m(\alpha,\sigma) = \sum_{i=1}^{m} c\left(\overline{y}_i, f_{\alpha,\overline{x},\sigma}(\overline{x}_i)\right) + \lambda\left\|f_{\alpha,x,\sigma}(\cdot)\right\|_{L_2}^2 =$$
$$= \sum_{i=1}^{m} c\left(\overline{y}_i, \sum_{j=1}^{m} \alpha_j k_{\sigma_j}(\overline{x}_i, \overline{x}_j)\right) + \lambda\alpha^T K\alpha \to \min_{\alpha,\sigma}, \quad (13)$$

To solve this problem one has to

1) Calculate the matrix $K = \left\{ \left\langle k_{\sigma_i}(\bar{x}_i, \cdot), k_{\sigma_j}(\bar{x}_j, \cdot) \right\rangle \right\}_{i,j=1}^m$ as a function of $\sigma = (\sigma_1, ..., \sigma_m)^{\mathrm{T}}$ ;

2) Solve the corresponding minimization problem (13) with respect to $\alpha = (\alpha_1, \alpha_2, ..., \alpha_m)^{\mathrm{T}} \in \mathbb{R}^m$ and $\sigma = (\sigma_1, \sigma_2, ..., \sigma_m)^{\mathrm{T}}$ .

Let us calculate matrix $K$ (9) in case of isotropic Gaussian kernel with variable width $\sigma$ :

$$k_\sigma(x, \bar{x}) = \frac{1}{(2\pi)^{n/2} \sigma^n} \exp\left\{ -\frac{\|x - \bar{x}\|^2}{2\sigma^2} \right\}.$$

Let

$$k_{\sigma_i}(x_i, x) = \frac{1}{(2\pi)^{n/2} \sigma_i^n} \exp\left\{ -\frac{\|x_i - x\|^2}{2\sigma_i^2} \right\}, \tag{14}$$

$$k_{\sigma_j}(x_j, x) = \frac{1}{(2\pi)^{n/2} \sigma_j^n} \exp\left\{ -\frac{\|x_j - x\|^2}{2\sigma_j^2} \right\}. \tag{15}$$

Then

$$\left\langle k_{\sigma_i}(x_i, \cdot), k_{\sigma_j}(x_j, \cdot) \right\rangle = \int_{\mathbb{R}^n} k_{\sigma_i}(x_i, x) k_{\sigma_j}(x_j, x) dx =$$

$$= \frac{1}{(2\pi)^n \sigma_i^n \sigma_j^n} \int_{\mathbb{R}^n} \exp\left\{ -\frac{\|x_i - x\|^2}{2\sigma_i^2} - \frac{\|x_j - x\|^2}{2\sigma_j^2} \right\} dx \tag{16}$$

$$= \frac{1}{(2\pi)^{n/2} \left( \sigma_j^2 + \sigma_i^2 \right)^{n/2}} \exp\left\{ -\frac{\|x_i - x_j\|^2}{\sigma_j^2 + \sigma_i^2} \right\}.$$

Proof. Indeed,

$$\frac{\|x_i - x\|^2}{2\sigma_i^2} + \frac{\|x_j - x\|^2}{2\sigma_j^2} = \frac{\sigma_j^2 \|x_i - x\|^2 + \sigma_i^2 \|x_j - x\|^2}{2\sigma_i^2 \sigma_j^2}$$

$$= \frac{\sigma_j^2 \left( \|x_i\|^2 + \|x\|^2 - 2\langle x_i, x \rangle \right) + \sigma_i^2 \left( \|x_j\|^2 + \|x\|^2 - 2\langle x_j, x \rangle \right)}{2\sigma_i^2 \sigma_j^2}$$

$$= \frac{\sigma_j^2 \|x_i\|^2 + \sigma_i^2 \|x_j\|^2 + \left( \sigma_j^2 + \sigma_i^2 \right) \|x\|^2 - 2\langle \sigma_j^2 x_i + \sigma_i^2 x_j, x \rangle}{2\sigma_i^2 \sigma_j^2}$$

$$= \frac{\sigma_j^2 \|x_i\|^2 + \sigma_i^2 \|x_j\|^2}{2\sigma_i^2 \sigma_j^2} + \frac{\sigma_j^2 + \sigma_i^2}{2\sigma_i^2 \sigma_j^2} \left( \|x\|^2 - 2\left\langle \frac{\sigma_j^2 x_i + \sigma_i^2 x_j}{\sigma_j^2 + \sigma_i^2}, x \right\rangle \right)$$

$$= \frac{\sigma_j^2 \|x_i\|^2 + \sigma_i^2 \|x_j\|^2}{2\sigma_i^2 \sigma_j^2} - \frac{\sigma_j^2 + \sigma_i^2}{2\sigma_i^2 \sigma_j^2} \left\| \frac{\sigma_j^2 x_i + \sigma_i^2 x_j}{\sigma_j^2 + \sigma_i^2} \right\|^2$$

$$+ \frac{\sigma_j^2 + \sigma_i^2}{2\sigma_i^2 \sigma_j^2} \left( x - \frac{\sigma_j^2 x_i + \sigma_i^2 x_j}{\sigma_j^2 + \sigma_i^2} \right)^2 ;$$

Now we simplify the expression

$$\frac{\sigma_j^2 \|x_i\|^2 + \sigma_i^2 \|x_j\|^2}{2\sigma_i^2 \sigma_j^2} - \frac{\sigma_j^2 + \sigma_i^2}{2\sigma_i^2 \sigma_j^2} \left\| \frac{\sigma_j^2 x_i + \sigma_i^2 x_j}{\sigma_j^2 + \sigma_i^2} \right\|^2 =$$

$$= \frac{\left( \sigma_j^2 \|x_i\|^2 + \sigma_i^2 \|x_j\|^2 \right)\left( \sigma_j^2 + \sigma_i^2 \right) - \left( \sigma_j^2 x_i + \sigma_i^2 x_j \right)^2}{2\sigma_i^2 \sigma_j^2 \left( \sigma_j^2 + \sigma_i^2 \right)} =$$

$$= \frac{\sigma_j^4 \|x_i\|^2 + \sigma_i^2 \sigma_j^2 \|x_j\|^2 + \sigma_i^2 \sigma_j^2 \|x_i\|^2 + \sigma_i^4 \|x_j\|^2 - \sigma_j^4 \|x_i\|^2 - \sigma_i^4 \|x_j\|^2 - 2\sigma_i^2 \sigma_j^2 \langle x_i, x_j \rangle}{2\sigma_i^2 \sigma_j^2 \left( \sigma_j^2 + \sigma_i^2 \right)} =$$

$$= \frac{+\sigma_i^2 \sigma_j^2 \|x_j\|^2 + \sigma_i^2 \sigma_j^2 \|x_i\|^2 - 2\sigma_i^2 \sigma_j^2 \langle x_i, x_j \rangle}{2\sigma_i^2 \sigma_j^2 \left( \sigma_j^2 + \sigma_i^2 \right)} = \frac{\|x_i - x_j\|^2}{\sigma_j^2 + \sigma_i^2} .$$

Finally, we obtain the required result:

$$\left\langle k_{\sigma_i}(x_i, \cdot), k_{\sigma_j}(x_j, \cdot) \right\rangle = \int_{\mathbb{R}^n} k_{\sigma_i}(x_i, x) k_{\sigma_j}(x_j, x) dx =$$

$$= \frac{1}{(2\pi)^n \sigma_i^n \sigma_j^n} \int_{\mathbb{R}^n} \exp \left\{ -\frac{\|x_i - x\|^2}{2\sigma_i^2} - \frac{\|x_j - x\|^2}{2\sigma_j^2} \right\} dx =$$

$$= \frac{1}{(2\pi)^n \sigma_i^n \sigma_j^n} \exp \left\{ -\frac{\|x_i - x_j\|^2}{\sigma_j^2 + \sigma_i^2} \right\} \int_{\mathbb{R}^n} \exp \left\{ -\frac{\sigma_j^2 + \sigma_i^2}{2\sigma_i^2 \sigma_j^2} \left( x - \frac{\sigma_j^2 x_i + \sigma_i^2 x_j}{\sigma_j^2 + \sigma_i^2} \right)^2 \right\} dx =$$

$$= \frac{1}{(2\pi)^n \sigma_i^n \sigma_j^n} \exp \left\{ -\frac{\|x_i - x_j\|^2}{\sigma_j^2 + \sigma_i^2} \right\} \int_{\mathbb{R}^n} \exp \left\{ -\frac{1}{2} \left( \frac{\sqrt{\sigma_j^2 + \sigma_i^2}}{\sigma_i \sigma_j} x - \frac{\sigma_j^2 x_i + \sigma_i^2 x_j}{\sigma_i \sigma_j \sqrt{\sigma_j^2 + \sigma_i^2}} \right)^2 \right\} dx =$$

$$= \frac{1}{(2\pi)^n \sigma_i^n \sigma_j^n} \exp \left\{ -\frac{\|x_i - x_j\|^2}{\sigma_j^2 + \sigma_i^2} \right\} \frac{\sigma_i^n \sigma_j^n}{\left( \sigma_j^2 + \sigma_i^2 \right)^{n/2}} \int_{\mathbb{R}^n} \exp \left\{ -\frac{1}{2} \left( x' - \frac{\sigma_j^2 x_i + \sigma_i^2 x_j}{\sigma_i \sigma_j \sqrt{\sigma_j^2 + \sigma_i^2}} \right)^2 \right\} dx' =$$

$$= \frac{1}{(2\pi)^{n/2} \left( \sigma_j^2 + \sigma_i^2 \right)^{n/2}} \exp \left\{ -\frac{\|x_i - x_j\|^2}{\sigma_j^2 + \sigma_i^2} \right\} .$$

**Remark 1.** Comparing problems (6) and (13), we can see that in the first case (6) the kernel matrix $K = \left\{ k(x_i, x_j) \right\}_{i,j=1}^{m}$ is given analytically but in the second case (13) we have to calculate it, i.e. to calculate inner products

$$\left\langle k_{\sigma_i}(x_i, \cdot), k_{\sigma_j}(x_j, \cdot) \right\rangle = \int_{\mathbb{R}^n} k_{\sigma_i}(x_i, x) k_{\sigma_j}(x_j, x) dx \, .$$

**Corollary 1.** For $n = 1$ formula (16) becomes

$$\left\langle k_{\sigma_i}(x_i, \cdot), k_{\sigma_j}(x_j, \cdot) \right\rangle = \frac{1}{(2\pi)^{1/2} \left( \sigma_i^2 + \sigma_j^2 \right)^{1/2}} \exp \left\{ -\frac{\left| x_i - x_j \right|^2}{2(\sigma_i^2 + \sigma_j^2)} \right\}, \quad \sigma_i > 0, \ \sigma_j > 0 \, .$$

Thus for an anisotropic kernel

$$k_{\sigma_i}(x, \overline{x}_i) = \frac{1}{(2\pi)^{n/2} \sigma_{i1} \times \sigma_{in}} \exp \left\{ -\sum_{l=1}^{n} \frac{\left| (x)_l - (\overline{x}_i)_l \right|^2}{2\sigma_{il}^2} \right\}, \quad \sigma_i = \left( \sigma_{i1}, \dots, \sigma_{in} \right),$$

it holds true

$$\left\langle k_{\sigma_i}(x_i, \cdot), k_{\sigma_j}(x_j, \cdot) \right\rangle = \frac{1}{(2\pi)^{n/2} \prod_{l=1}^{n} \left( \sigma_{il}^2 + \sigma_{jl}^2 \right)^{1/2}} \exp \left\{ -\sum_{l=1}^{n} \frac{\left| (x_i)_l - (x_j)_l \right|^2}{\sigma_{il}^2 + \sigma_{jl}^2} \right\}.$$

**Remark 2.** Let $\rho(x) \geq 0$, $x \in \mathbb{R}$, be some probabilistic density function, $\int_{\mathbb{R}} \rho(x) dx = 1$. Consider a family of kernels $k_\sigma(x, \overline{x}) = \sigma^{-1} \rho \left( \sigma^{-1} (x - \overline{x}) \right)$, $\sigma > 0$, and inner products

$$\left\langle k_{\sigma_1}(x, \overline{x}_1), k_{\sigma_2}(x, \overline{x}_2) \right\rangle = \int_{\mathbb{R}} k_{\sigma_1}(x, \overline{x}_1), k_{\sigma_2}(x, \overline{x}_2) dx = \varphi(\overline{x}_1, \overline{x}_2, \sigma_1, \sigma_2) \, .$$

Function $\varphi(\overline{x}_1, \overline{x}_2, \sigma_1, \sigma_2)$ may not be known analytically but can be calculated (tabulated) in advance and then can be used for quick construction of the large kernel matrix (9).

### 3.3 Analytical solution of the quadratic regression problem

Analogously to [17, 29], we can obtain a closed form solution of problem (11):

$$\alpha = \left( K + \lambda I \right)^{-1} \overline{y} \, , \tag{17}$$

where matrix $K = \left\{ k_{ij} \right\}_{i,j=1}^{m} = \left\{ \left\langle k_{\sigma_i}(\overline{x}_i, \cdot), k_{\sigma_j}(\overline{x}_j, \cdot) \right\rangle \right\}_{i,j=1}^{m}$ is symmetric and positive definite, $I$ is the diagonal identity matrix, $\overline{y} = \left( \overline{y}_1, \dots, \overline{y}_m \right)^{\mathrm{T}}$. The optimal value is $F_{\min}^{m} = \lambda \overline{y}^{\mathrm{T}} \left( K + \lambda I \right)^{-1} \overline{y}$. In case of kernels (14), (15), coefficients $k_{ij}$ are given by formula (16).

# 4 Adaptation of kernels to the problem data

Adaptation of kernels to problem data here means solution of problem (12) over kernel parameter $\sigma$. In this paper we illustrate this on the Example of Section 2.

## 4.1 Adaptation of kernels in case of Nadaraya-Watson regression

As noted in subsection 2.1, NW-regression is rather robust to the choice of kernel width, since it just averages values $\{\bar{y}_i\}_{i=1}^m$ in this or that way. The kernel width influences on the number of essential summands in (2). Below, on Figures 3, 4 we illustrate NW-regression (3) with $k = 10$ and $k = 3$. There is no essential difference between obtained approximations.

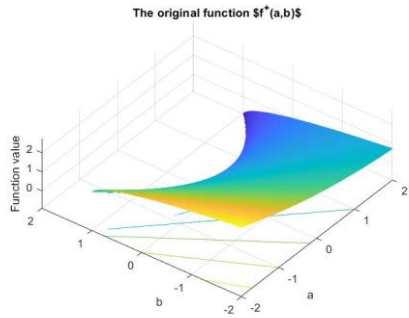

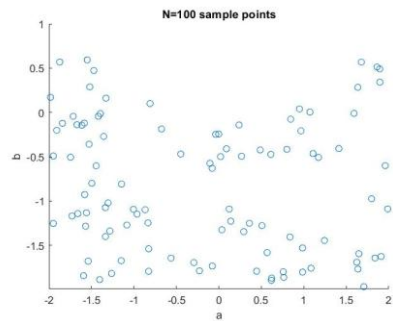

Fig. 1. The original function $f^+(a,b)$ points within feasible set $(a,b) \in [-2,2]^2$ :

Fig. 2. $m = 100$ random sample within the feasible set.

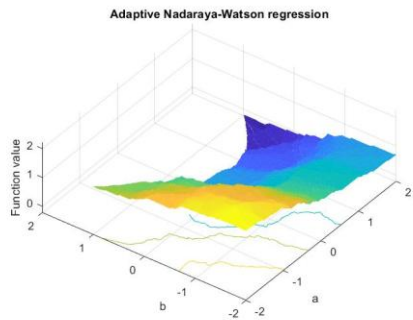

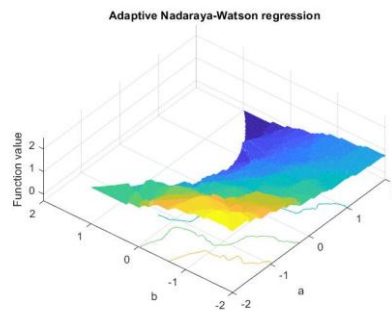

Fig. 3. NW-regression with $k = 10$ summands and absolute accuracy $\Delta^m \approx 0.7162$ .

Fig. 4. NW-regression with $k = 3$ summands and absolute accuracy $\Delta^m \approx 0.6688$

## 4.2 Adaptation of kernels in case of kernel SVM

In case of kernel SVM the kernel width is quite essential for the quality of approximation of the original regression function. Partially, this drawback is coped with the choice of regularization parameter $\lambda$ in (10), which prevents over fitting of the regression model. However, it is still necessary to coordinate kernel width with data, at least with the number $m$ of used observations $\{(\overline{x}_i, \overline{y}_i)\}_{i=1}^m$. A too small kernel width leads to oscillating regression function but a too large kernel width leads to loosing local information in data. So the possible solution is to find optimal widths by solving problem (13) on additional sample of data. So we can use additional data not only for selection of the regularization level but also for choosing kernels. Remark that works [19, 22] use the same initial data set $\{(\overline{x}_i, \overline{y}_i)\}_{i=1}^m$ for finding optimal common width through leave-one-out of $\{(\overline{x}_i, \overline{y}_i)\}_{i=1}^m$ cross validation technique.

The key problem in the kernel SVM method is choosing the value of the regularization parameter $\lambda$. For this we use the cross validation technique, i.e. minimization over $\lambda$ of the approximation error on an additional test data. We illustrate this scheme on Example of Section 2.

### 4.2.1. Separate adaptation of kernel weights $(\alpha_1, ..., \alpha_m)$ and a common kernel width $\sigma$

For a numerical illustration consider the following special case of problem (18):

$$\sum_{i=1}^m \left( \overline{y}_i - \sum_{j=1}^m \frac{\alpha_j}{(2\pi)^{n/2} \sigma^n} \exp\left\{ -\frac{\left\| \overline{x}_j - \overline{x}_i \right\|^2}{2\sigma^2} \right\} \right)^2 +$$

$$+ \sum_{i,j=1}^m \frac{\alpha_i \alpha_j}{(2\pi)^{n/2} \left(2\sigma^2\right)^{n/2}} \exp\left\{ -\frac{\left\| \overline{x}_j - \overline{x}_i \right\|^2}{2\sigma^2} \right\} \to \min_{\alpha_1, ..., \alpha_m; \sigma > 0}.$$

In this case we first solve problem (11) and find optimal weights $\alpha^0 = \left\{ \alpha_i^0 \right\}_{i=1}^m$ under some fixed initial value $\sigma = \sigma^0 > 0$, then for the found $\alpha^0$ and a new sample $\left\{ (x_i', y_i') \right\}_{i=1}^{m'}$ solve problem (12) and find under fixed $\alpha = \alpha^0$ a common optimal width $\sigma^1 > 0$. Finally, for the found $\sigma^1$ and the original sample $\{(\overline{x}_i, \overline{y}_i)\}_{i=1}^m$ again solve problem (11) and find new weights $\alpha^1 = \left\{ \alpha_i^1 \right\}_{i=1}^m$, and so on. Now the kernel SVM regression function takes on the form $f^m(x) = \sum_{i=1}^m \frac{\alpha_i^s}{(2\pi)^{n/2} (\sigma^s)^n} e^{-0.5 \|x - x_i\|^2 / (\sigma^s)^2}$, $s \geq 1$. Here the optimal weights $\alpha_i^s$ and optimal width $\sigma^s$ depend on the regularization parameter $\lambda$. Denote $\Delta^m = \sup_{a,b} \left| f^+(a,b) - f^m(a,b) \right|$ the accuracy estimate. We illustrate this approach on the Example. The following figures 5-7 illustrate the proposed procedure of calculating $\sigma^1, \alpha^1$ for $\lambda = 10^{-4}$. Figure 5 shows that the regularized losses have a unique minimum at some $\sigma^1$. This observation is the main experimental finding of this section. Figures 6, 7 display the obtained function estimates for some initial $\sigma^0$ and for the optimal width $\sigma^1$.

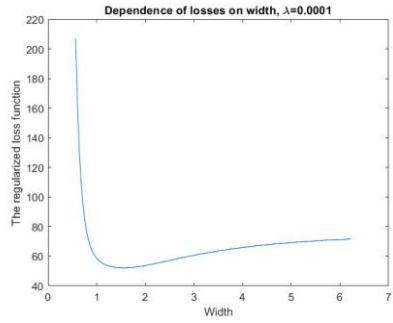

Fig. 5. Dependence of the losses on
the kernel width for an additional sample.

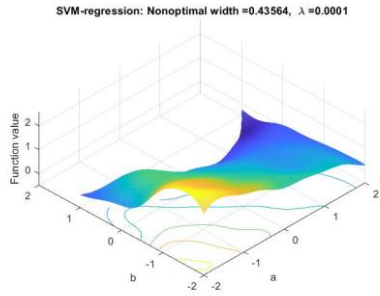 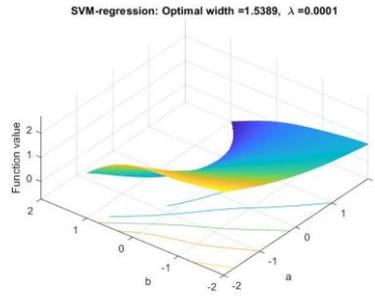

Fig. 6. SVM-regression for $\lambda = 0.001$
and initial width $\sigma^0 = 0.4356$
with absolute accuracy $\Delta^m = 1.4593$.

Fig. 7. SVM-regression for $\lambda = 10^{-4}$
and optimized width $\sigma^1 = 1.5389$
with absolute accuracy $\Delta^m = 0.2275$.

The next figures show the dependence of accuracy, optimal width on the regularization parameter $\lambda$. It can be seen that the minimal approximation error in Fig. 8 is achieved at $\lambda = 10^{-4}$, for which width achieves its maximum value in Fig. 9.

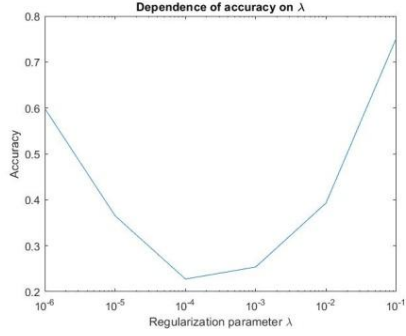 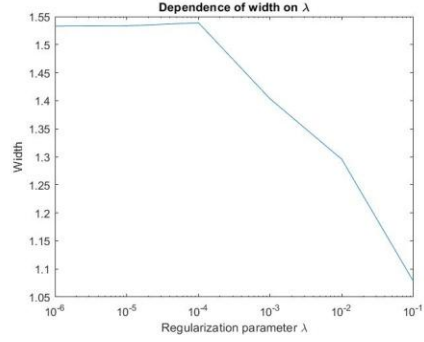

Fig. 8. Dependence of accuracy on $\lambda$ .    Fig. 9. Dependence of width on $\lambda$ .

### 4.2.2. Adaptation of kernel weights $\left(\alpha_1,...,\alpha_m\right)$ and widths $\left(\sigma_1,...,\sigma_m\right)$

In this case we solve problem (10), which now has the form

$$\sum_{i=1}^{m}\left(\overline{y}_i - \sum_{j=1}^{m}\frac{\alpha_j}{(2\pi)^{n/2}\sigma_j^n}\exp\left\{-\frac{\left\|\overline{x}_j - \overline{x}_i\right\|^2}{2\sigma_j^2}\right\}\right)^2 +$$

$$+\lambda\sum_{i,j=1}^{m}\frac{\alpha_i\alpha_j}{(2\pi)^{n/2}\left(\sigma_i^2+\sigma_j^2\right)^{n/2}}\exp\left\{-\frac{\left\|\overline{x}_j - \overline{x}_i\right\|^2}{2\left(\sigma_i^2+\sigma_j^2\right)}\right\}\rightarrow\min_{\alpha_1,...,\alpha_m;\sigma_1>0,...,\sigma_m>0}$$

through solving problems (11) and (12) in turn iteratively several times. At iteration $s$ we solve problem (11) under fixed $\sigma=\sigma^{s-1}=\left(\sigma_1^{s-1},...,\sigma_m^{s-1}\right)^{\mathrm{T}}$,

$$R^m(\alpha^s,\sigma^{s-1})=\min_{\alpha}\left[R^m(\alpha,\sigma^{s-1})=\left(\sum_{i=1}^{m}c\left(\overline{y}_i,f_{\alpha,\overline{x},\sigma^{s-1}}(\overline{x}_i)\right)+\lambda\Psi\left(f_{\alpha,\overline{x},\sigma^{s-1}}(\cdot)\right)\right)\right]$$

to find optimal weights $\alpha^s=\left\{\alpha_i^s\right\}_{i=1}^{m}$, and the corresponding function $f_{\alpha^s,\overline{x},\sigma^{s-1},\lambda}(\cdot)$, $s\geq 1$. Then for the found fixed $\alpha=\alpha^s$ and a new sample $\left\{(x_i',y_i')\right\}_{i=m}^{m'}$, $m'>m$, we solve problem (12),

$$R^m(\alpha^s,\sigma^s)=\min_{\sigma}\left[R^m(\alpha^s,\sigma)=\left(\sum_{i=1}^{m'}c\left(\overline{y}_i',f_{\alpha^s,\overline{x},\sigma}(\overline{x}_i')\right)+\lambda\Psi\left(f_{\alpha^s,\overline{x},\sigma}(\cdot)\right)\right)\right],$$

and find optimal widths $\sigma^s = \left(\sigma_1^s, ..., \sigma_m^s\right)^{\mathrm{T}}$, $s \geq 1$, and function $f_{\alpha^s, \overline{x}, \sigma^s, \lambda}(\cdot)$; and so on. Remark that in case of the quadratic error function $c(y - f) = (y - f)^2$ and the regularization $\Psi(f) = \sum_{i,\,j=1}^{m} \alpha_i \alpha_j k_{ij}$, problem (11) admits analytical solution (17). More over, experiments show that in (12) we can use the same (training) data set $\left\{(x_i, y_i)\right\}_{i=1}^{m}$ as in (11). The test data set can be used for selection of optimal regularization parameter $\lambda$ by minimizing the test error over $\lambda$ (for some $s$),

$$\Delta^s(\lambda_{opt}) = \inf_{\lambda > 0} \max_{x \in \left\{x'_j\right\}_{j=m}^{m'}} \left| f(x) - f_{\alpha^s, \overline{x}, \sigma^s, \lambda}(x) \right|.$$

We illustrate the proposed adaptation procedure on Example from Section 2 for $m = 50,\ 100$ random training sample points $\left\{x_i\right\}_{i=1}^{m}$ and random test points $\left\{x_i\right\}_{i=m+1}^{2m}$ in the feasible set $\left\{x = (a, b) : a^2 - 4b \geq 0\right\}$. The initial width $\sigma = \sigma_i^0$ of the Gaussian kernel $k_\sigma(x - x_i)$ was taken as 4-th nearest to $x_i$ point from $\left\{x_i\right\}_{i=1}^{m}$. Figures 10, 11 show the accuracy of obtained kernel regression function $f_{\alpha^s, \overline{x}, \sigma^s, \lambda}(\cdot)$ after $s = 15$ iterations on training points, test points, and the actual accuracy for different $\lambda \in [0.0001, 0.1]$. Figures 12, 13 show function $f_{\alpha^{15}, \overline{x}, \sigma^0, \lambda_{opt}}(\cdot)$, i.e. without adaptation of widths, and $f_{\alpha^{15}, \overline{x}, \sigma^{15}, \lambda_{opt}}(\cdot)$, i.e., with adaptation of widths, $\lambda_{opt} = 0.005$.

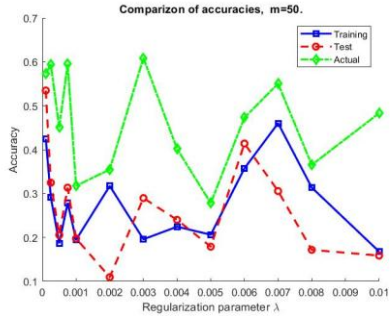

Fig. 10.

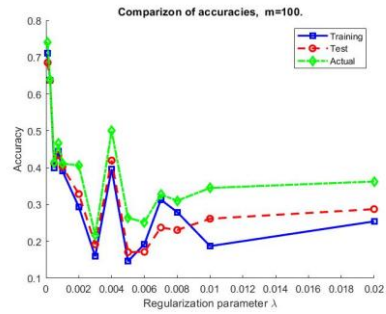

Fig. 11.

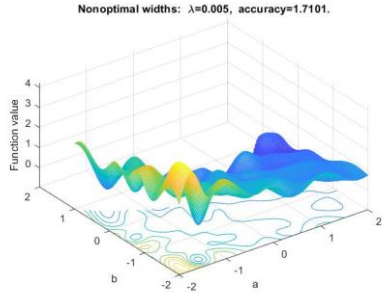

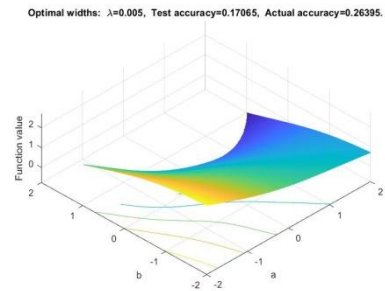

Fig. 12. Regression function before adaptation of widths.

Fig. 13. Regression function after adaptation of widths.

Thus the main experimental observation of this subsection is that joint minimization of the regularized error function indeed improves the accuracy of the obtained kernel SVM regression function.

## 5 Conclusions

Any large scale applied mathematical models contain numerous parameters, see e.g., [34-40]. A standard approach to analyze such models is to investigate the dependences (robustness) of results of modeling on parameters. However these dependences may be very complex, nonconvex and nonsmooth. One possible approach to reconstruct these dependences is to use the kernel regression method presented in the present paper. To improve the standard kernel regression method we extend it by adapting kernels, e.g. their widths, to the sample data. We illustrate the approach on Nadaraya-Watson and SVM regressions. In the present paper it was checked a limited adaptation, namely in the case of Nadaraya-Watson regression we used several nearest-neighbor data points. It appears that the accuracy weakly depends on the number accounted neighbor data points. In the case of SVM regression we used additional minimization of the regularized error function over multiple kernel widths. And even in a limited adaptation of widths, this trick considerably improves the accuracy of the SVM method.